\documentclass[12pt]{amsart} \usepackage{amssymb}

\usepackage{helvet} 

\usepackage{courier} 

\reversemarginpar 
 \normalmarginpar 

\let\oldmarginpar\marginpar 
\renewcommand\marginpar[1]{\-\oldmarginpar{\raggedright\small\sf #1}}

\title{Congruences for rational points on  varieties over finite fields}

\author{N.~Fakhruddin and C.~S.~Rajan}
 
\address{School of Mathematics, Tata Institue of Fundamental Research,
Homi Bhabha Road, Mumbai 400005, India}

\email{naf@math.tifr.res.in, \;rajan@math.tifr.res.in}

\newcommand{\nc}{\newcommand}

 \newcommand{\Q}{\mathbb{Q}}
 \newcommand{\C}{\mathbb{C}}
 \newcommand{\A}{\mathbb{A}}
 \nc{\W}{\mathbb{W}}
\newcommand{\mc}{\mathcal} \newcommand{\mb}{\mathbb}

\nc{\G}{\Gamma}

\renewcommand{\P}{\mb{P}}

\nc{\aff}{{\A}^1} \nc{\naive}{\!\sim_n} \nc{\Spec}{\mathrm{Spec}}

\nc{\wt}{\widetilde} \nc{\wh}{\widehat}

\nc{\sr}{\stackrel}

\nc{\lr}{\longrightarrow} \nc{\wb}{\overline}
\newtheorem{thm}{Theorem}[section] \newtheorem{prop}[thm]{Proposition}
 \newtheorem{cor}[thm]{Corollary}
\newtheorem{lem}[thm]{Lemma}

\theoremstyle{definition} \newtheorem{defn}[thm]{Definition}
\newtheorem{ques}[thm]{Question} \newtheorem{rem}[thm]{Remark}
 \input{xypic}
\begin{document}

\begin{abstract}
We show that the number of rational points on the fibres of a proper
morphism of smooth varieties over a finite field  $k$ whose generic
fibre has a ``trival'' Chow group  of zero cycles is congruent to $1
\mod |k|$. As a consequence we  prove that there is a rational point
on any degeneration of  a smooth proper rationally chain connected
variety over a  finite field. We also obtain a generalisation of the
Chevalley-Warning theorem. 
\end{abstract}

\maketitle

\section{Introduction}

\subsection{}
The  main result of this paper is the following:
\begin{thm}
Let $f_i: X_i \to Y$, $i = 1,2$ be  proper dominant
morphisms of smooth irreducible varieties over a  finite field $k$ and
let $g: X_1 \to X_2$ be a dominant morphism over $Y$.  Let $Z_i$ be
the generic fibre of $f_i$ and assume that
$g_*:CH_0({Z_1}_{\wb{k(X_1)}})_{\Q}\to CH_0({Z_2}_{\wb{k(X_1)}})_{\Q}$
is an isomorphism.  Then for all $y \in Y(k) $,
\[
|f_1^{-1}(y)(k)| \equiv |f_2^{-1}(y)(k)| \mod |k| .
\]
\label{thm:MT}
\end{thm}

Specialising to the case $X_2=Y$ and $f_2 = Id_Y$,
we obtain:
\begin{cor}
Let $f: X \to Y$ be a proper dominant morphism of smooth irreducible
varieties over a  finite field $k$. Let $Z$ be the generic fibre of
$f$ and assume that $CH_0(Z_{\wb{k(X)}})_{\Q} = \Q$.  Then for all $y
\in Y(k) $, $|f^{-1}(y)(k)| \equiv 1 \mod |k|$.
\label{cor:CW}
\end{cor}

When $Y = Spec(k)$, the above corollary reduces to a theorem
of H.~Esnault \cite{esnault-fano}. Since we do not assume that $f$
is smooth, we are able to obtain congruences even for singular
varieties.

\smallskip

An immediate consequence of Corollary \ref{cor:CW} is that  for $f:X \to
Y$ as above, $|X(k)| \equiv |Y(k)| \mod |k|$. By taking $Y$ to be
a point in Theorem \ref{thm:MT} we obtain: 

\begin{cor}
Let $g:X_1 \to X_2$ be a dominant morphism of smooth proper varieties over
a finite field $k$. If $g_*: CH_0(X_{\wb{k(X)}})_{\Q} \to
CH_0(Y_{\wb{k(X)}})_{\Q}$ is an isomorphism, then $|X(k)| \equiv
|Y(k)| \mod |k|$.\footnote{This can also be proved using the method of \cite{kahn-cong}.}
\label{cor:CW2}
\end{cor}

Simple examples show that in contrast to Corollary \ref{cor:CW}, the
assumption that $X_1$ and $X_2$ be proper cannot be omitted here.

\smallskip



The triviality of the Chow group of zero cycles of
degree $0$, or even rational chain connectedness\footnote{We learnt of
 such an example, due to  J.~Koll\'ar, from a talk by J.~Iyer at the
 University of Essen in May 2003; but see Remark \ref{rem:chain}}
is not sufficient to
guarantee the existence of a rational point for proper varieties over
finite fields which are not smooth.  However, using alterations, we
give a criterion for the existence of rational points  which can be
applied to all degenerations of smooth rationally chain connected
varieties.

\begin{cor}
Let $f:X \to Y$ be a proper dominant morphism of
irreducible varieties  over a finite field $k$ with  $Y$  smooth. Let
$Z$ be the generic fibre of $f$ and assume that
$Z$ is smooth and $CH_0(Z_{\wb{k(X)}})_{\Q} = \Q$.  Then for any $y \in Y(k)$,
$f^{-1}(y)(k) \neq \emptyset$.
\label{cor:nonempty}
\end{cor}


The corollary below generalises the Chevalley--Warning theorem
 \cite{greenberg}, which is the special case  $P = \P^n$ and $L_i =
 \mc{O}(d_i)$ with $\sum d_i \leq n$.
\begin{cor}
Let $P$ be a smooth projective geometrically irreducible variety over
a finite field $k$. Let $ L_1, \cdots, L_r$ be  very ample line
bundles on $P$ such that
$(K_{P}\otimes L_1\otimes \cdots\otimes L_r)^{-1}$ is ample,
 where $K_P$ is the canonical bundle of
$P$. For $i = 1,2, \dots, r$, let  $s_i \in H^0(P,L_i)$. Then
\[
\big |\big \{x \in P(k) \big | s_i(x) = 0 \ , \ i=1,2, \dots, r \big \} \big |
\equiv 1 \mod  |k| \ .
\]
\label{cor:CW3}
\end{cor}

Note that the congruence formula of  Katz \cite{sga72},  when it
applies, only gives congruences modulo $p = char(k)$.
It would be interesting to know whether an analogue of the Ax--Katz
theorem \cite{katz-ax} holds in the above situation or whether all low
degree intersections as above over $C_1$ fields always have rational
points.


\bigskip

Our proofs, as is Esnault's, are based on the method of decomposing
the diagonal originated by Bloch, the Grothendieck-Lefschetz trace formula
 and  Deligne's
integrality theorem. 
The novelty of our approach lies in the use of proper
correspondences which allows us to deal with non-proper $Y$; this is
crucial  for obtaining congruences for (or even just existence of)
rational points on singular varieties. The main technical ingredient
for this is the refined intersection theory of Fulton and MacPherson
\cite{fulton-it}.

\section{Proper Correspondences}

\subsection{}
In this section we introduce  the concept of proper correspondences 
and prove some of its properties.
\footnote{A.~Nair has informed us that a variant of this definition
has been  considered, in a different context, by E.~Urban in the
preprint: Sur les repr\'esentations $p$-adiques associ\'ees aux
repr\'esentations cuspidales de $GSp_4/\Q$, 2001.}
  Most of the
proofs are simple modifications of those in Fulton's book
\cite{fulton-it} and we prove only what we need for later use; several
other results in \cite[Chapter 16]{fulton-it} for (usual)
correspondences have analogues for proper  correspondences.

\begin{defn}
Let $X$ and $Y$ be smooth irreducible varieties over a field $k$.  The
group of \emph{proper correspondences} from $X$ to $Y$, $PCorr(X,Y)$,
is the free abelian group generated by irreducible subvarieties
$\Gamma \subset X \times Y$ which are proper over $Y$ modulo the
subgroup generated by elements of the form
\[
\{div(f)| f \in k(Z)^*, \ Z \subset X \times Y \mbox{ irreducible and
proper over }Y\} \ .
\]
\end{defn}

\smallskip
This is clearly a graded abelian group; we shall grade it by dimension
(lower indices) or codimension (upper indices) 
as is convenient. For a  cycle $\Gamma$ (not
necessarily irreducible) in the free abelian group as above we  shall
denote its class in $PCorr(X,Y)$ by $[\G]$.

If $f:X \to Y$ is a proper morphism, then the graph of $f$ gives an
element $[\Gamma_f]$ of $PCorr(X,Y)$. We shall show below that proper
correspondences induce maps on cohomology generalising the maps
induced by proper morphisms.

\subsection{}
We recall some properties of \'etale (co)homology and cycle class maps
that we shall need. The main reference is \cite[Expos\'es
VI-IX]{douady-verdier}; we note that the quasi-projectivity hypotheses
there can be removed using the methods of  \cite{fulton-it}. The
compatibility of refined intersection products and refined cycle class
maps stated below can be deduced using the methods of \cite[Chapter
19]{fulton-it}; see also \cite{corti-hanamura}.

\bigskip
Let $K$ be an algebraically closed field and fix a prime number $l \neq
char(K)$.  For a variety $X$ over $K$ let $H^i(X):=H^i_{et}(X, \Q_l)$,
$H^i_c(X):= H^i_{c, et}(X, \Q_l)$ and let $H^i_{Z}(X):=H^i_{Z, et}(X,
\Q_l)$ for $Z$  a closed subvariety of $X$.  We also let $H_i(X)$
denote the locally finite (or Borel-Moore) homology $H^{et}_i(X,
\Q_l)$. For any of the above groups $H$, we denote by $H(n)$, for an integer $n$,
the corresponding Tate twisted group.


For a cycle $\alpha = \sum_i a_i Z_i$,  on a variety $X$, we denote by
$|\alpha|$  the support $\cup_i Z_i$ of $\alpha$.  The following
properties are proved in the references cited above:
\begin{enumerate}
\item (Projection formula) For any variety $X$ there are cap product
maps  $H^i(X)\otimes H_j(X) \to H_{j-i}(X)$ such that  if $f: X \to Y$
is a proper morphism, $u \in H^i(Y), v \in H_j(X)$ then $f_*(f^*(u)
\cap v) = u \cap f_*(v)$ in $H_{j-i}(Y)$.
\item For $Z$ an irreducible variety of dimension $n$, $H_{2n}(Z)(-n)$ is
one dimensional.
\item For $\alpha$  a cycle of dimension $k$ on a variety $X$  the
  fundamental class defines a
canonical element $cl(\alpha) \in H_{2k}(|\alpha|)(-k)$. This  maps
to an element, also denoted by $cl(\alpha)$, in $H_{2k}(X)(-k)$ which is
zero if $\alpha = 0 \in CH_k(X)$.

\item  For $X$  smooth of pure dimension
$n$  there is a canonical isomorphism $H_{2k}(|\alpha|)(-k) \cong H^{2n -
2k}_{|\alpha|}(X)(n-k)$, so we also get an element $cl(\alpha) \in H^{2n -
2k}_{|\alpha|}(X)(n-k)$.

\item For $f: X \to Y$  a morphism and $\alpha$  a cycle on $X$ 
of dimension $n$ whose
support is proper over $Y$, $f_*(cl(\alpha)) = cl(f_*(\alpha))$ in
$H^*_{f(|\alpha|)}(Y)(-n)$.

\item For $\alpha$ and $\beta$  cycles of pure dimension $k$ and $l$
respectively in  a smooth irreducible variety $X$ of dimension $n$
\[ cl(\alpha \cdot \beta) = cl(\alpha) \cup cl(\beta) 
\in H^{4n-2(k+l)}_{|\alpha| \cap |\beta|}(X)(n-k-l)\ ,
\]
 where the product on the left is the refined intersection product.
\end{enumerate}

\subsection{}

Let  $\dim(X) = n, ~\dim(Y) = m$ and let 
$\Gamma$ be a representative for an element of $PCorr_r(X,Y)$. 
Then $\Gamma$ induces a linear map
$\Gamma_*:H^i(X) \to H^{2m -2r +i}(Y)(m-r)$ as the composite of the
sequence of linear maps:
\begin{multline*}
H^i(X) \stackrel{p_X^*}{\lr} H^i(X \times Y) \stackrel{\cdot
cl(\Gamma)}{\lr} H^{i + 2(n +m -r)}_{|\Gamma|}(X \times Y) \cong \\
H_{2r
-i}(|\Gamma|)(-r) 
\stackrel{{p_Y}_*}{\lr} H_{2r -i}(Y)(-r) \cong H^{2m -2r +
i}(Y)(m-r) \ .
\end{multline*}

Since $X$ and $Y$ are smooth, by duality 
we also get maps $\Gamma^*: H^i_c(Y) \to H^{2n -2r +
i}_c(X)(n-r)$.

Now suppose that $k$ is a finite  field and $K$  an algebraic closure
of $k$.   All the cohomology
groups discussed above, for varieties over $K$ which are base changed
from varieties over $k$, have a continuous action of $Gal(K/k)$. The
maps $\Gamma^*$ and $\Gamma_*$ are then  compatible with the Galois action. 

\smallskip
The following 
lemma implies that the maps $\Gamma_*$ and $\Gamma^*$ only depends
on $[\Gamma] \in PCorr(X,Y)$.
\begin{lem}
Let $\Gamma$ be as above. Suppose there exists a closed subvariety $Z
\subset X \times Y$ such that $|\Gamma| \subset Z$, $[\Gamma] = 0$ in
$CH_*(Z)$ and $p_Y|_Z:Z \to Y$ is proper. Then $\Gamma^*$ and
$\Gamma_*$ are the zero maps.
\label{lem:rateq}
\end{lem}
\begin{proof}
Suppose $Z$ is a closed subvariety of $X\times Y$ such that the
projection $p_Y|_Z:Z \to Y$ is proper. It follows from the projection
formula that $\G_*$ can also be defined as, 
\[
H^i(X) \stackrel{p_X|_{Z}^*}{\lr} H^i({Z}) \stackrel{\cdot
cl(\Gamma)}{\lr}  H_{2r
-i}(Z)(-r) \stackrel{{p_Y}_*}{\lr} H_{2r -i}(Y)(-r) \cong H^{2m -2r +
i}(Y) (m-r),
\]
where  $cl(\Gamma)$ is considered as an element in
 $H_{2r}(Z)(-r)$. The lemma follows from the fact that $cl(\Gamma) = 0$ in
 $H_{2r}(Z)(-r)$.
The statement for $\G^*$ follows by duality.
\end{proof}

\subsection{}
Let $X,Y,Z$ be smooth irreducible varieties over $K$ and let
$[\Gamma_1]$ (resp. $[\Gamma_2]$) be a proper correspondence from $X$
to $Y$ (resp. $Y$ to $Z$).  Analogous to the definition of
composition of correspondences \cite[Chapter 16]{fulton-it},  we define
$[\G_2] \circ [\G_1] \in PCorr(X,Z)$ by
\[ 
[\G_2] \circ [\G_1] = [{p_{XZ}}_*({p_{XY}}^*(\G_1) \cdot
{p_{YZ}}^*(\G_2))] \ ,
\]
 where the $p$'s denote the projection morphisms from $X \times Y
\times Z$ and the product is the refined intersection product. Note
that ${p_{XY}}^*(\G_1) \cdot {p_{YZ}}^*(\G_2)$ is a cycle which is
well defined upto rational equivalence on $|\G_1| \times_Y |\G_2|$.
Since $|\G_1| \times_Y |\G_2|$ is proper over $Z$, its image in $X
\times Z$ is also proper over $Z$, so $[\G_2] \circ [\G_1]$ is a well
defined element of $PCorr(X,Z)$.

\begin{lem}
Let $X,Y,Z$ and $\G_1,\G_2$ be as above. Then $(\G_2 \circ \G_1)_* =
(\G_2)_* \circ (\G_1)_*$ as maps from $H^*(X)$ to $H^*(Z)$ and $(\G_2
\circ \G_1)^* = (\G_1)^* \circ (\G_2)^*$ as maps from $H^*_c(Z)$ to
$H^*_c(X)$.
\label{lem:funct}
\end{lem}
 
\begin{proof}
The key point is the compatibility of the refined cycle class maps
with refined intersection products.

Let $a \in H^*(X)$. Then
\begin{align*}
(\G_2)_* \circ (\G_1)_*(a) & = {p_Z^{YZ}}_*(cl(\G_2) \cdot
{p_Y^{YZ}}^*({p_Y^{XY}}_*(cl(\G_1) \cdot{p_X^{XY}}^*(a)))) \\ &=
{p_Z^{YZ}}_*(cl(\G_2) \cdot {p_{YZ}^{XYZ}}_*({p_{XY}^{XYZ}}^*(cl(\G_1)
\cdot{p_X^{XY}}^*(a))))\\ &=
{p_Z^{YZ}}_*({p_{YZ}^{XYZ}}_*({p_{YZ}^{XYZ}}^*(cl(\G_2)) \cdot
{p_{XY}^{XYZ}}^*(cl(\G_1) \cdot{p_X^{XY}}^*(a))))\\ &=
{p_Z^{YZ}}_*({p_{YZ}^{XYZ}}_*({p_{YZ}^{XYZ}}^*(cl(\G_2)) \cdot
{p_{XY}^{XYZ}}^*(cl(\G_1)) \cdot{p_X^{XYZ}}^*(a))))\\ &=
{p_Z^{XZ}}_*({p_{XZ}^{XYZ}}_*({p_{YZ}^{XYZ}}^*(cl(\G_2)) \cdot
{p_{XY}^{XYZ}}^*(cl(\G_1)) \cdot{p_X^{XYZ}}^*(a))))\\ &=
{p_Z^{XZ}}_*(cl(\G_2 \circ \G_1) \cdot{p_X^{XZ}}^*(a))\\ &= (\G_2
\circ \G_1)_*(a)
\end{align*}
We use the projection formula several times along with compatibility
of the cycle class map with smooth pullbacks, products and proper
pushforwards.
\end{proof}

\subsection{}

The key technical result which allows us to deduce congruences from
cycle theoretic information is the following:

\begin{prop}
Let $X,Y, \G$ be as above and assume that $m =n = r$. If
$\dim(p_X(|\G|)) < n$ then all the eigenvalues of (the geometric)
Frobenius acting on $\G^*(H^i_c(Y)) \subset H^i_c(X)$ are algebraic
integers which are divisible by $|k|$.
\label{lem:div}
\end{prop}

\begin{proof}
Replacing $k$ by a finite extension does not affect the conclusion of
the lemma, so without loss of generality we may assume that $\G$ is a
geometrically irreducible subvariety of $X \times Y$. Using the
definition of $\G_*$ as given in the proof of Lemma \ref{lem:rateq},
$\G_*$ is
the composite of the following maps:
\[
H^i(X) \sr{p_X^*}{\lr} H^i(\G) {\lr} H_{2n-i}(\G)(-n)
\sr{{p_Y|\G}_*}{\lr} H_{2n-i}(Y)(-n) \cong H^{i}(Y) \ , 
\]
where we have used the hypothesis that $m=n=r$. 

Suppose $\pi:\G' \to \G$ is a proper dominant  generically finite
morphism with $\G'$ smooth and geometrically irreducible. Then the
projection formula shows that we may replace $\G$ with $\G'$ and $p_X$
(resp. $p_Y$) with $p_X \pi$ (resp. $p_Y \pi$) in the above without
changing the image of the composite.

Let $W$ be the Zariski closure of $p_X(\G)$ in $X$ and let $\dim(W) =
t$.  By the theorem of De Jong \cite[4.1]{dejong-alt1} we may find
$\G'$ as above, $\sigma:W' \to W$ with $\sigma$ proper dominant
generically finite and $W'$ smooth, and $p: \G' \to W'$ such that $p_X
\pi = i_W \sigma p$, where $i_W:W \to X$ is the inclusion.  Since $p_X \pi =
(i_W \sigma) p$ and both $\G'$ and $W'$ are smooth, it follows from
the functoriality of pushforward maps that $({p_X} \pi)_* = (i_W
\sigma)_* p_*$. By
the remarks of the previous paragraph we see that the
image of $\G_*$ is the same as that of the composite of the sequence:
\[
H^i(X) {\lr} H^i(W') {\lr} H^i(\G')  {\lr}
H^i(Y) \ .
\]
Since $\G'$ and $W'$ are smooth, by duality we get a factorisation of
$\G^*$ as a composite of maps, 
\[
H^j_c(Y) {\lr} H^j_c(\G') {\lr} H^{2(t-n)+j}_c(W')(t-n) {\lr} H^j_c(X) \ .
\]
By Deligne's integrality theorem \cite{sga72}, Expos\'e XXI,
 Corollaire 5.5.3, all the eigenvalues of Frobenius on $H^*_c(W')$ are
 algebraic integers. Since $n-t > 0$ and the geometric Frobenius acts
 on $\Q_l(t-n)$ by $|k|^{n-t}$,  the proposition follows.
\end{proof}

\section{Proofs of the main results}
Using the results of the previous section, we now give the proofs of
the results stated in the introduction. 
\subsection{}

\begin{proof}[Proof of Theorem \ref{thm:MT}]
If $Y$ is not geometrically irreducible, then $Y(k) = \emptyset$ so
there is nothing to prove. The hypothesis on the Chow group  implies
that the geometric connected components of $Z_1$ and $Z_2$ are in
bijective correspondence. Since $X_1$ and $X_2$ are smooth, we can 
 assume that $X_1$, $X_2$  and $Y$ are  geometrically
irreducible.

Let $W$ be  an irreducible  subvariety of $X_1$
which maps generically finitely and dominantly to $X_2$ and let
$d$ be the degree of $W$ over $X_1$. Let $\G_g$ be
the graph of $g$ in $X_1 \times_k X_2$, let $\G_W$ be the transpose of the
graph of $g|_W$ embedded in $X_2\times_k X_1$ and let $\Gamma_1 = (\G_W \circ
\G_g)/d$. Since $W$ is a subvariety of $X_1$, $\G_W$ is a proper
correspondence, and $\G_g$ is a  proper correspondence since $g$ is
proper. Hence $\G_1$ is a proper correspondence of dimension
$n_1=\dim{X_1}$. 

Let $V_2$ be the open subset of $X_2$ over which
$g|_{W}$ is finite and let $V_1 = g^{-1}(V_2)$. By the construction of $W$,
$p_{12}^*(\Gamma_g)$ and
$p_{23}^*(\Gamma_W)$, which are subvarieties of
$X_1 \times_k X_2 \times_k X_1$, meet properly when pulled back to $V_1 \times_k X_2 \times_k X_1$.
It follows that one can write $\Gamma_1 = \G_1' + \G_1''$ in $PCorr(X_1,X_1)$
where $(Id_{X_1} \times g)_*(\G_1') = \G_g$ in $PCorr(X_1,X_2)$ and
$p_1(|\G_1''|)$ is a proper subvariety of $X_1$.


Let $\Delta_{X_1}$ be the diagonal in $X_1 \times_k X_1$
and consider the proper correspondence  
$\Gamma_2:=\Delta_{X_1}-\Gamma_1$ and its restriction (i.e.~flat pullback)
$\gamma_2$ to $Z_1 \times_{k(Y)} Z_1$.
Since $(Id_{X_1} \times g)_*(\Delta_{X_1}) = \G_g$,
 it follows from the previous paragraph that $(Id_{Z_1} \times g|_{Z_1})_*(\gamma_2)$
can be represented by a  cycle on $Z_1 \times_{k(Y)} Z_2$ which becomes zero
when restricted to $k(Z_1) \times_{k(Y)} Z_2$. Since the map 
$g_*:CH_0({Z_1}_{\wb{k(X_1)}})_{\Q}\to CH_0({Z_2}_{\wb{k(X_1)}})_{\Q}$
is injective, it follows that $\gamma_2$ can be represented by a cycle
on $Z_1 \times_{k(Y)} Z_1$ whose support maps to a proper subvariety of
$Z_1$ by the projection to the first factor.
By taking the Zariski closures in $X_1\times_Y X_2$, it follows that 
$\Gamma_2$ can be represented by a proper correspondence
on $X_1 \times_k X_1$ whose support maps to a proper subvariety of $X_1$
by the projection to the first factor.

By Proposition \ref{lem:div} all the eigenvalues of Frobenius
acting on the image of $\Gamma_2^*$ in $H^*_c(X_1)$ are divisible by
$|k|$. It follows from Lemma \ref{lem:funct} and the definition of $\Gamma_1$,
 that the image
of $\Gamma_1^*$ is contained in the image of $g^*$.
By the definition of $\Gamma_2$, we conclude that
all the eigenvalues of Frobenius
acting on the cokernel of $g^*:H^*_c(X_2) \to H^*_c(X_1)$ are divisible by $|k|$.

The hypotheses of the theorem, and the above discussion, are not
affected if we replace $Y$ by an open subvariety $U$ and $X_i$ by 
$f_i^{-1}(U)$, $i=1,2$, so we may  assume that $Y$ has only one rational point
$y$.  Applying the Grothendieck-Lefschetz trace formula and the
statement on eigenvalues above, we  conclude that the $|X_1(k)| \equiv |X_2(k)|
\mod |k|$. Each rational point of $X_i$ must lie over the unique
rational point $y$ of $Y$,  hence the theorem follows. 

\end{proof}

\begin{cor}
Let $f: X \to Y$ be a proper dominant generically finite  morphism of
 smooth irreducible varieties over a  finite field $k$ such that the
 extension of function fields $k(Y) \to k(X)$ is purely inseparable.
 Then for any $y \in Y(k) $, $|f^{-1}(y)(k)| \equiv 1 \mod
 |k|$. \label{cor:birat}
\end{cor}

\begin{proof}
Since $(Z_{\wb{k(X)}})_{red}$ is isomorphic to $Spec(\wb{k(X)})$, the
hypothesis on $CH_0$ is trivially satisfied.
\end{proof}

\begin{rem}
Since we only use intersection theory (resp.~cohomology groups) with
rational (resp.~$\Q_l$) coefficients, the above results hold even
when $X$ and $Y$ are quotients of smooth varieties by finite
groups. \label{rem:1}
\end{rem}




\subsection{}
\begin{proof}[Proof of Corollary \ref{cor:nonempty}]
As in the proof of Theorem \ref{cor:CW}, we may assume that $X$ and
$Y$ are geometrically irreducible and that $Y$ has a unique rational
point.

By a result of De Jong \cite[5.15]{dejong-alt2} there exists an
irreducible  variety $X'$ over $k$ which is the quotient of a smooth
variety by a finite group and a proper dominant generically finite
morphism $\pi:X' \to X$ such that the extension of function fields
$k(X) \to k(X')$ is purely inseparable. Let $f' = f \pi :X' \to Y$ and
let $Z'$ be the  generic fibre of $f'$.
Since $X'$ is irreducible over $k$, $Z'$ is irreducible over
$k(Y)$. $Z$ is geometrically irreducible since it is smooth over
$k(Y)$; since the extension of function fields induced by the map $Z'
\to Z$ is $k(X) \to k(X')$, $Z'$ is also geometrically irreducible.
The induced morphism $({Z'}_{\wb{k(X)}})_{red} \to Z_{\wb{k(X)}}$ thus
satisfies the assumptions of Lemma \ref{lem:insep} below,  so
$CH_0({Z'}_{\wb{k(X')}} )_{\Q} = \Q$.  Applying Corollary \ref{cor:CW}
(cf. Remark \ref{rem:1}) to the morphism $f':X' \to Y$ we see that
$X'(k) \neq \emptyset$. Thus $X(k) \neq \emptyset$.
\end{proof}

\begin{rem}
For singular $X$ it is not  always true that $|f^{-1}(y)(k)| \equiv 1
\mod |k|$, however the only examples we know where this fails are
non-normal.
\end{rem}

\begin{rem}
\label{rem:chain}
$X(k) \neq \emptyset$ if $X$ is a proper variety over a finite field
$k$ which is rationally connected i.e.  any two general points of
$X(\Omega)$ are contained in an irreducible rational curve in $X$
defined over $\Omega$, where $\Omega \supset k$ is a universal domain.
(To prove this we may replace $X$ by $X'$ as in the proof of Corollary
\ref{cor:nonempty}. Since the extension of function fields $k(X) \to
k(X')$ is  purely inseparable, it follows that $X'$ is also rationally
connected, so  $CH_0({X'}_{\wb{k(X')}})_{\Q} = \Q$. By Remark
\ref{rem:1} it follows that $X'(k) \equiv 1 \mod |k|$, hence $X'(k)
\neq \emptyset$. Thus $X(k) \neq \emptyset$.) However, a degeneration
of a rationally connected variety is in general only rationally chain
connected, so one cannot use this to prove Corollary
\ref{cor:nonempty} even if $Z$ is rationally connected.
\end{rem}

\begin{lem}
\label{lem:insep}
Let $f:X \to Y$ be a proper dominant morphism of irreducible varieties
over  an algebraically closed field $K$. Assume that $Y$ is
smooth,  $f$ is generically finite and
the extension of function fields $K(Y) \to K(X)$ is purely
inseparable.  Then $f_*:CH_0(X)_{\Q} \to CH_0(Y)_{\Q}$ is an
isomorphism.
\end{lem}

\begin{proof}
$f_*$ is surjective because $f$ is surjective. Using the refined
intersection theory of \cite{fulton-it} and the hypothesis on $Y$, we
see that there is a natural map $f^*:CH_0(Y)_{\Q} \to CH_0(X)_{\Q}$.
By the assumptions on $f$ there exists an open subset $U$ of $Y$
such that for $V = f^{-1}(U)$, $f|_V:V \to U$ is a bijection. By
the moving lemma (which is trivial for zero cycles) $CH_0(X)$
(resp. $CH_0(Y)$) is generated by the closed points in $V$
(resp. $U$). This shows that 
$f^*$ is a surjection since for $y
\in U$,  $f^*([y]) = [f^{-1}(y)]$.
Since $f_*f^*$ is multiplication by $\deg(f)$, it
follows that $f_*$ is an isomorphism.
\end{proof}


\subsection{}
\begin{proof}[Proof of Corollary \ref{cor:CW}]
Let $Y = \prod_i H^0(P, L_i)$ and let $X \subset Y \times
P$ be the zero scheme of the universal section of $\bigoplus_i
p_P^*L_i$.   Since the $L_i$ are basepoint
free, the map $X \to P$ is smooth, hence $X$ is also smooth. Since the
$L_i$ are assumed to be very ample, Bertini's theorem
implies that the generic fibre $Z$ of the projection $f:X \to Y$ is
smooth. The assumptions on the $L_i$ and the adjunction formula
imply that $Z$ is a Fano variety. 
 By a theorem of Campana
\cite{campana-fano} and Kollar-Miyaoka-Mori \cite{kmm}  $Z$ is
rationally chain connected, so  $CH_0(Z_{\wb{k(X)}})_{\Q} = \Q$.  The
proof is concluded by applying Corollary \ref{cor:CW} to $f$.
\end{proof}

\begin{rem}
The proof shows that instead of assuming that the $L_i$
are very ample it suffices to assume that they are basepoint free and
that $Z$ is smooth.
\end{rem}

A motivic version of the above corollary in the case of hypersurfaces
in $\P^n$ is proved in \cite{bel}.


\section{Further questions}
\label{sec:questions}
It seems likely that the following mixed characteristic analogue of
Corollary \ref{cor:nonempty} has a positive answer:
\begin{ques}
let $K$ be a finite extension of $\Q_p$, $\mc{O}$ its ring of
integers, $k$ the residue field and $X$  a smooth proper variety over
$K$  such that $CH_0(X_{\wb{K(X)}})_{\Q} = \Q$.  If $\mc{X} \to
Spec(\mc{O})$ is a proper scheme with generic fibre isomorphic to $X$
and closed fibre $X_0$, then is $X_0(k) \neq \emptyset$?
\label{ques:nonempty}
\end{ques}

From our proofs we do not obtain any
information about the valuations of the eigenvalues of Frobenius
acting on the \'etale cohomology of the fibres of $f$. One is thus
lead to ask the following:
\begin{ques}
Let $f: X \to Y$ be a proper dominant morphism of smooth irreducible
varieties over a  finite field $k$. Let $Z$ be the generic fibre of
$f$ and assume that $CH_0(Z_{\wb{k(X)}})_{\Q} = \Q$.  Then for all $y
\in Y(k) $ and $i>0$, does $|k|$ divide all the eigenvalues of
Frobenius acting on $H^i((f^{-1}(y))_{\wb{k}}, \Q_l)$ ($l \neq char(k)$)?
\label{ques:div}
\end{ques}

Theorem \ref{thm:MT}
has the following Hodge theoretic analogue:
\begin{thm}
Let $f_i: X_i \to Y$, $i=1,2$, be  proper dominant morphisms of smooth irreducible
varieties over  $\C$ and let $g:X_1 \to X_2$
be a dominant morphism over $Y$.  Let $Z_i$ be the generic fibre of $f_i$ and assume
that $g_*: CH_0({Z_1}_{\wb{\C(X_1)}})_{\Q} \to CH_0({Z_2}_{\wb{\C(X_1)}})_{\Q}$
is an isomorphism. Then $gr^F_0(Coker(g^*:
H^i_c(X_1,\C) \to H^i_c(X_2,\C))) = 0$ for all $i$, where
$F^{\overset{\mbox{\huge{.}}}{ }}$ denotes the Hodge filtration of
Deligne.
\label{thm:hodge}
\end{thm}
The proof  is essentially the same as that of Theorem \ref{thm:MT}  so
we omit the details: one only needs to replace Lemma \ref{lem:div} by
its Hodge theoretic counterpart.

\smallskip
The  Hodge version of Question \ref{ques:div} is:
\begin{ques}
Let $f: X \to Y$ be a proper dominant morphism of smooth irreducible
varieties over  $\C$.  Let $Z$ be the generic fibre of $f$ and assume
that $CH_0(Z_{\wb{\C(X)}})_{\Q} = \Q$.  Then for all $y \in Y(\C) $
and $i>0$, is $gr^F_0(H^i(f^{-1}(y),\C)) = 0$?
\label{ques:hodge}
\end{ques}


It seems likely, as was suggested to us by P.~Brosnan, that there
should be a purely  motivic  statement from which Theorems
\ref{cor:CW} and \ref{thm:hodge} should follow after taking
realisations. The category of motives over a base of Corti and
Hanamura \cite{corti-hanamura} would seem to be a natural choice in
which to formulate such a statement.

\begin{rem} After this  paper was submitted both Question \ref{ques:div} and
(a stronger form of) Question \ref{ques:hodge} 
were shown to have positive answers by H.~Esnault; her results appear in
the appendix to this article \cite{esnault-app}.
 A positive answer to Question \ref{ques:nonempty}
in case $\mc{X}$ is regular was also obtained by her in \cite{esnault-deligne}; the results
of that paper include a stronger version of Corollary \ref{cor:CW} in the case that
$Y$ is a curve.
\end{rem}


{\bf{\emph{Acknowledgements.}}}  N.F.~thanks Arvind Nair for a
remark which convinced him that Corollary \ref{cor:CW} should be
true (before Theorem \ref{thm:MT} was proved)
and Patrick Brosnan, H\'el\`ene Esnault, Madhav Nori and V.~Srinivas  for
their comments on the results of this paper. He also thanks
H\'el\`ene Esnault and Marc Levine for an invitation to visit the
University of  Essen in May--June 2003 where he was supported by the
DFG and via the Wolfgang Paul prize program of the Humbolt Stiftung;
lectures of H\'el\`ene Esnault that he attended at that time provided
part of the motivation for the results and questions in this paper.

We also thank 
Bruno Kahn for his comments on the first version of
this paper and for suggesting Theorem \ref{thm:MT} as a common
generalisation  of Corollaries \ref{cor:CW} and \ref{cor:CW2}.

\bigskip

\end{document}